\documentclass[11pt,a4paper]{article}

\usepackage{amssymb}
\usepackage{amssymb,amsmath} %,amscd,concmath}
\usepackage{comment}
\usepackage{graphicx}
\usepackage{subfig}
\usepackage{epsfig}
\usepackage{cite}
\usepackage{multicol}
\usepackage{color}
\usepackage{framed}
\usepackage{textcomp}
\usepackage{stmaryrd} % for \varotimes
\usepackage{hyperref}

\usepackage[hmargin=2.3cm,vmargin=3.5cm]{geometry}
\usepackage{enumerate}
\usepackage{bm} % for bold math symbols ?
\usepackage{amsmath}
\usepackage{amsthm} % for a prove environment
\usepackage{amstext}
\usepackage{mathrsfs} % for script math symbols: \mathscr{...}
\usepackage{pifont} % for \hand symbol
\usepackage{url} % for urls in the bibliography
\newcommand{\kommentar}[1]{}

\newcommand{\brackets}[1]{\left ( #1 \right )}

\newcommand{\converg}{\rightarrow}

\newcommand{\dtilde}[1]{\tilde{\tilde{#1}}}
\newcommand{\NVert}{\mathcal{N}}

\newcommand{\FT}{\mathcal{F}}

\newcommand{\Real}{\mathbb{R}}

\newcommand{\Complex}{\mathbb{C}}

\newcommand{\Nat}{\mathbb{N}}
\newcommand{\Whole}{\mathbb{Z}}

\theoremstyle{definition}

\newcommand{\gausslower}[1]{\left\lfloor #1 \right\rfloor}

\newcommand{\const}{\text{const}}

\newcommand{\abs}[1]{\left \lvert#1 \right \rvert}

\newcommand{\norm}[1]{\lVert#1\rVert}

\newcommand{\E}{\text{E}}

\renewcommand{\d}{\,\text{d}}

\newcommand{\ONotSym}{\mathcal{O}}

\renewcommand{\d}{\,\text{d}}

\newcommand{\keywords}[1]{\textbf{Keywords:} #1}

\newcommand{\AMSclassification}[1]{\textbf{AMS 2010 Subject Classification:} #1}

\newcommand{\multisetZ}{\{Z\}}

\newcommand{\conv}{*}
\newcommand{\SchwartzSp}{\mathcal{S}}
\newcommand{\DiracComb}{\raisebox{-0.2mm}{\text{$\mathsf{L\!L\!I}$}}}

\newcommand{\si}{\text{si}}

\excludecomment{sectionremoved}
\excludecomment{textremoved}
\excludecomment{mycommenttext} % this defines an environment whose content later is regarded as comment, so do not define mycommenttext later as env.

\newcounter{TheoremCounter}
\newcounter{LemmaCounter}
\newcounter{PropositionCounter}
\newcounter{RemarkCounter}
\newenvironment{mytheorem}{\refstepcounter{TheoremCounter}\mbox{ }\\\textbf{Theorem \arabic{TheoremCounter}}:}{}
\newenvironment{mylemma}{\refstepcounter{LemmaCounter}\mbox{ }\\\textbf{Lemma \arabic{LemmaCounter}}:}{}

\newenvironment{myremark}{\refstepcounter{RemarkCounter}\mbox{ }\\\textbf{Remark \arabic{RemarkCounter}}:}{}

\newcounter{AlgorithmCounter}

\makeatletter
\renewcommand\section{\@startsection
  {section}{1}{0mm}%name, level, indent
  {-\baselineskip}%             beforeskip
  {0.3\baselineskip}%            afterskip
  {\normalfont\Large\bfseries}}% style
\renewcommand\subsection{\@startsection
  {subsection}{2}{0mm}%name, level, indent
  {-\baselineskip}%             beforeskip
  {0.1\baselineskip}%            afterskip
  {\normalfont\normalsize\bfseries}}% style
\makeatother

\begin{document}

\title{Efficient computation of the cumulative distribution function of a linear mixture of independent random variables}
\author{Thomas Pitschel\footnote{Correspondence address: th (dot) pitschel (at) tu-braunschweig (dot) de}}
\maketitle

\begin{abstract}
For a variant of the algorithm in \cite{Pitschel2019DeterministicBootstrapping} 
(arxiv.org/abs/1903.10816) to compute the approximate density
or distribution function of a linear mixture of independent random variables known by a 
finite sample, it is presented a proof of the functional correctness,
i.e. the convergence of the computed distribution function towards the true distribution
function (given the observations) as the algorithm resolution is 
increased to infinity. 
The algorithm (like its predecessor version) bears elements which 
are closely related to early known methods for numerical inversion
of the characteristic function of a probability distribution, however
here efficiently computes the \emph{complete} distribution function.
Possible applications are in computing 
the distribution of the bootstrap estimate in any linear bootstrap method (e.g. 
in the block bootstrap for the mean as parameter of interest, or residual bootstrap in linear regression with fixed design),
or in elementary analysis-of-variance hypothesis testing.
\end{abstract}

\mbox{ }\\
\mbox{ }\\
\AMSclassification{62G07, 62G09} 

\section{Introduction}
\label{sec:Introduction}
Bootstrapping is a resampling technique employed to estimate a
parameter $\theta$ of a distribution (most often variances, in order to systematically
construct confidence intervals) in the face of no a-priori 
knowledge about the true distribution from which the available data
points $X_i$, $i=1\dots n$, (in this text real-valued) are deemed to be originating 
from. In particular, bootstrapping is capable of extending hypothesis testing beyond the 
requirement that a parametric form of the distribution of the involved noise
is known. 

The conventional way of using bootstrapping in practice is to evaluate
a certain estimator $T_n^*$ on multiple ''synthetic'' data sets (=replicates)
which are generated (in the basic type) by randomly sampling data points from
the orginal data set, with the aim of letting them mimick this original
data set with regards to its distributional characteristics. The
so derived values will, under suitable assumptions placed on the bootstrap
method (i.e. on estimator $T_n^*$ and the resampling distribution) as well
as on the parameter of interest, mimick the distribution of an estimator
$T_n$ which is a centered and scaled derivate of the estimator $\hat\theta$
of the actually targetted parameter $\theta$. (=''consistency'') The method 
thus allows, under suitable assumptions, to find the distribution \mbox{of $\hat\theta$}
(when applied to the whole population) and thus the ultimately targetted confidence 
regions.

In \cite{Pitschel2019DeterministicBootstrapping}, an algorithm was outlined which allows to circumvent the
above described ''random element'' in conventional bootstrap method usage and arrive at an approximation of the 
distribution of $T_n^*$ (conditional on the available data points $X_i$) by directly (deterministically) concluding from the $X_i$,
provided that the bootstrap method used is linear in the sense 
defined in \cite{Pitschel2019DeterministicBootstrapping}.
In the present text, the aim is to prove rigorous statements on the
error of the approximation, and the focus will be on an in $N$ asymptotic 
result. The main statement on this is contained in section \ref{sec:AsymptoticsBasedOnTildefRepresentation}, as well as in the appendix.

Besides facilitating practical realization of the bootstrapping 
procedure for linear bootstrap methods, the computation of the distribution of a linear 
mixture of random variables potentially has applications related to other 
estimation problems. 
In the wider context of parameter estimation in econometric models at small
sample sizes, \cite{Phillips1982ExactSmallSampleTheoryDiscPaper621} 
surveys methods for estimating the distribution of estimators in
simultaneous equations models.
\cite{Davies1980TheDistributionOfALinearCombination} has error bounds 
on the (for this field relevant) computation of the distribution of 
quadratic forms of multiple independent normally distributed random variables,
which may be written as linear combinations of independent $\chi^2$ variables, i.e.
which are in a form amenable to the here examined algorithm.
In \cite{AbateChoudhuryWhitt1999AnIntroductionToNumericalInversionAndApplications},
methods are developed for inverting transforms of probability distributions, and
applied to obtain tail probabilities in queuing systems. (Compared to their work,
which is geared towards obtaining single probability values with high accuracy and typically uses
analytical expressions of the characteristic functions of the underlying pdfs, the 
algorithm here computes the whole distribution function (sampled equidistantly).)

The field of most direct application could initially appear to be the hypothesis 
testing for the comparison of two populations' means, %(each of which a finite sample is observed), 
which occurs throughout social and life sciences, see for example \cite{MunzelHauschke2003NonparametricTestForProvingNoninferiority}, 
\cite{HarkinWebb2016MonitoringGoalProgress}. However, the recognized established
non-parametric tests (Mann-Whitney test, resp. rank-based Kruskal-Wallis test in the more-than-two-groups case) 
will be preferable due to the more rigorous statements derivable at equal assumptions.
Such test scenarios may be regarded as analysis-of-variance with discretely levelled factors.
For confidence intervals in factor models with continuously valued factors, 
applicability of the here examined algorithm remains an open question.

The text proceeds as follows: In the next section, the problem to be tackled is
restated, and elementary definitions and relations given. In section \ref{sec:AsymptoticResult}, 
elementary properties relevant to the inversion of the characteristic functions are stated. In section \ref{sec:FiniteNresults}, 
an alternative explicit expression of the algorithm output is derived
and the main convergence result proven with it.

\section{Setting and common definitions}
\label{sec:CommonDefinitions}
Let $X_i$, $i=1\dots n$, be real-valued observations deemed to be realizations from
some (not further considered) random variable, and let $\hat F_n$ be the associated empirical distribution function,
i.e. $\hat F_n(x) = n^{-1} \cdot \sum_{i=1}^n 1_{X_i \leq x}$. Let the random variable $X$ be defined as
distributed according to $\hat F_n$. Let $X^{[j]}$, $j=1\dots m$, be independent
random variables distributed as $X$. Let $a_j \in \Real$, $j=1\dots m$.
Define 
\begin{align}
    Z := \sum_{j=1}^m a_j X^{[j]}  \label{eqn:Z_definition}
\end{align}
and denote by $F_Z$ its cumulative distribution function. In \cite{Pitschel2019DeterministicBootstrapping}, 
it was stated that the algorithm presented there, here called Algorithm 1,
computes an approximation to the density $f_Z$ of $Z$. With $G$ and $g_k$ defined as\footnote{It is ''$i$'' after ''$2\pi$'' or in denominator the imaginary unit, otherwise usually acts as index variable.}
\begin{align}
    G_{a_j X^{[j]}}(\nu) & = \frac{1}{n} \sum_{i=1}^n \exp(- 2\pi i \cdot a_j X_i \nu) \label{eqn:G_ajXj_definition} \\
    G(\nu) & = \prod_{j=1}^m G_{a_j X^{[j]}}(\nu) \label{eqn:G_definition}, \quad \quad g_k = G(k \Delta \nu),        % \label{eqn:g_k_definition}
\end{align}
i.e. $g_k$ as computed in Algorithm 1, set for $x \in \Real$ and $N \in \Nat$
\begin{align}
    \tilde h(x) & := \frac{1}{N} \sum_{k=0}^{N-1} g_k \exp(2 \pi i \cdot x k/T)    \label{eqn:tilde_h_x_definition} \\
    \tilde f(x) & := \frac{1}{N} \sum_{k=-N+1}^{N-1} g_k \exp(2 \pi i \cdot x k/T)   \label{eqn:tilde_f_x_definition} \\
    \hat f_Z(x) & := \tfrac{N}{T} \cdot \tilde f(x).
\end{align}
It is $g_0=1$ and $\abs{g_k} \leq 1$; because of $g_0=1$ and $g_{-k} = g_k^*$, one easily finds %\\[1mm]
$\tilde f(x) = 2 \cdot \text{Re}(\tilde h(x)) - 1/N$.
For an $i \in \{0, \dots, N-1\}$, it is $\tilde f_i := \tilde f(i \tfrac{T}{N})$ (called $f_i$ in equation
(5) in \cite{Pitschel2019DeterministicBootstrapping}).
Set $\tilde h_i := \tilde h(i \tfrac{T}{N})$. One comfortably proves 
$\sum_{i=0}^{N-1} \tilde h_i = 1$ and therefore also $\sum_{i=0}^{N-1} \tilde f_i = 1$.

An alternative representation of $F_Z$ can be stated by recognizing
that $Z$ takes finitely many values. Denote by $\{Z\}$ the set of those values,
and let $p_z$ the probability mass of value $z \in \{Z\}$ (as implied by $\hat F_n$). Then
\begin{align}
    F_Z = \sum_{z \in \{Z\}} p_z \cdot 1_{z \leq x}.
\end{align}

\section{Towards convergence proof: Elementary properties of $G$} % "$N$-asymptotics, route1"
\label{sec:AsymptoticResult}
In this section, the proof of the convergence of the algorithm result to 
the desired distribution $F_Z$ is prepared by recollecting 
some elementary properties for the inversion of a characteristic 
function and stating relevant definitions.
''Convergence'' and ''asymptotic'' here refers to 
the behaviour as $N \converg \infty$, while (on the contrary) the $n$, i.e. the number of 
data points, remains fixed. Insofar, it is examined here the necessary 
computational resources to be expended to achieve a sufficiently accurate 
result, ideally independent of the given fixed input sample size $n$.

In prospect of wanting to apply the knowledge of the distribution to
deriving confidence intervals and rejection probabilities, it is the error
in estimating $F_Z(z)$ for each $z \in \Real$ that is of interest.

Let $G(\nu)$ be as in eqn. \eqref{eqn:G_definition}. Remark \ref{lem:G_is_like_characteristic_function} below asserts
that $G$ equals --up to argument-side scaling-- the characteristic function
of the distribution of $Z$. Since the characteristic function determines the 
distribution uniquely, it is $G(\nu)$ representing this distribution exactly.
Applying then an inversion formula, the earliest version of which appears to have 
been stated by L{\'e}vy and a modified version of which is stated below as Lemma \ref{lem:InversionFormulaForG},
therefore yields the true cumulative distribution function $F_Z$ of $Z$.
    \kommentar{
    Then, from Lemma 1 and the
    existence of an inverse map (from characteristic function to distribution function),
    one concludes that $G(\nu)$ is characterizing the distribution of $Z$ exactly.
    Applying then this inverse map to $G$, such the one given by Levy's inversion 
    formula \cite{Levy1925CalculDesProbabilites} (below stated as Lemma 2)
    therefore yields the true cumulative distribution function $F_Z$ of $Z$.
    }

\begin{myremark}
\label{lem:G_is_like_characteristic_function}
Let $X_i, i=1\dots n,$ be real-valued observations, and $\hat F_n$
the associated empirical distribution function. Let random variable $Z$ be as in equation \eqref{eqn:Z_definition}, again 
with the contained random variables $X^{[j]}$ independent and distributed according to $\hat F_n$.
Let $G$ be defined as in \eqref{eqn:G_definition}.
Then the characteristic function of $Z$ (conditional on $X_1, \dots, X_n$), defined as 
$t \mapsto \E(\exp(i t Z))$, % $| \{X_i\})$, 
is given by $t \mapsto G(- t/(2\pi))$.
\end{myremark}

The statement of the remark is based on the well-known argument that by 
independence of the $X^{[j]}$, the expectation separates into a product
of expectations. Since the $X^{[j]}$ are discretely distributed, each of the 
expectations can be written as in \eqref{eqn:G_ajXj_definition}. \qed

\mbox{ }\\
The following lemma is stated here for completeness, but will not be made use of directly. %It will not be made use of.
The expressions derive from the well-known inversion formulae for the characteristic 
function originally appearing in \cite{Levy1925CalculDesProbabilites}, %(L{\'e}vy 1925), 
by simple variable substitution and
noting that $F_Z(x_0) = 0$.

\begin{mylemma}
\label{lem:InversionFormulaForG}
Let $X_i$, $\hat F_n$, $Z$ and $G$ as in Lemma 1. Let $x_0 < \min(Z)$.
Then at every point $x$ of continuity of $F_Z$, it is $F_Z(x) = \bar F(x)$ with
\begin{align}
    \bar F(x) := \frac{1}{2\pi} \int_{-\infty}^\infty \frac{e^{2 \pi i \cdot x \nu} - e^{2 \pi i \cdot x_0 \nu}}{i \nu} G(\nu) \d \nu, \label{eqn:bar_F_x_definition}
\end{align}
and because $G(-\nu) = G^*(\nu)$, at the same points also $F_Z(x) = 2 \cdot \text{Re}(H(x))$ with
\begin{align}
     H(x) & := \frac{1}{2 \pi} \int_{0}^\infty \frac{e^{2 \pi i \cdot x \nu} - e^{2 \pi i \cdot x_0 \nu}}{i \nu} G(\nu) \d \nu  \label{eqn:H_x_definition} \\
     F_Z(x) = \bar F(x) & = \frac{1}{\pi} \int_{0}^\infty \frac{\text{Im}\brackets{(e^{2 \pi i \cdot x \nu} - e^{2 \pi i \cdot x_0 \nu}) \cdot G(\nu)}}{\nu}  \d \nu. \label{eqn:F_Z_integral_repres}
\end{align}
Here, since $\abs{G(\cdot)}$ is bounded and the fraction in the integrand of \eqref{eqn:H_x_definition}
tends to $2 \pi (x-x_0)$ as $\nu \converg 0$, the integrand is continuous and bounded everywhere and 
the integral in \eqref{eqn:H_x_definition} can well be evaluated as Riemann integral around $\nu=0$. 
\end{mylemma}
\qed

\begin{mycommenttext}
The inversion formula is quite standard. A form of the inversion formula closer to the one 
used here is found also in \cite{Billingsley1979ProbabilityAndMeasure} (Billingsley 1979) or in \cite{Phillips1982ExactSmallSampleTheoryDiscPaper621} (Phillips 1982).
\end{mycommenttext}

\begin{mycommenttext}
Note: For random variables with existing probability density function, the
characteristic function essentially equals the Fourier transformation of that
density. Correspondingly, the density is reconstructed by applying the inverse Fourier 
transform on the characteristic function. It follows that the cumulative distribution
function must be obtainable by integrating the inverse Fourier transform.
\end{mycommenttext}

One route of analysis at this point would proceed by using Lemma \ref{lem:InversionFormulaForG}
to establish a link between an estimate of the distribution function (defined in the second next paragraph)
and the true $F_Z$ (as given by its integral representation above). For the purpose of this 
text however, a more accessible and possibly more illuminating route is preferred.
 
The following definitions are used throughout the remainder of the text.
By $\tilde F_i$ (see equation \eqref{eqn:tilde_F_i_definition} below) 
it is denoted the ''cyclical'' accumulative sum 
of the $\tilde f_i$ of Algorithm 1. %(which are called ''$f_i$'' in \cite{Pitschel2019DeterministicBootstrapping}, equation (5) there). 
It maps \mbox{''$\!\!\!\!\mod N$''} into the integer range $\{0, \dots, N-1\}$.

Set $T_Z := \max Z - \min Z$. Let $T > T_Z$, $\kappa := T / T_Z$, assume $z_{min} < 0 < z_{max}$ 
and let $i_{min} = \gausslower{\tfrac{N}{T} \cdot \kappa \cdot z_{min}}$ and then $x_0 = i_{min} T/N$.
For %$i=0 \dots N-1$ and 
$i' = i_{min} \dots (i_{min}+N-1)$ set
\begin{align}    
    \tilde F_{i'} & := \sum_{i''=i_{min}}^{i'-1} \tilde f_{i'' \text{mod} N}  = \frac{1}{N} \sum_{k=-N+1}^{N-1} G(k \Delta \nu) \sum_{i''=i_{min}}^{i'-1} e^{2 \pi i \cdot \tfrac{i'' k}{N}}     \label{eqn:tilde_F_i_definition} \\
                & = \frac{1}{N} \sum_{k=-N+1}^{N-1} G(k \Delta \nu) \sum_{i''=0}^{i'-i_{min}-1} \exp(2 \pi i \cdot \tfrac{(i_{min} + i'') k}{N})   \label{eqn:tilde_F_i_evaluation}
\end{align}

\section{Results at finite $N$ and asymptotic result}
\label{sec:FiniteNresults}
The route of analysis pursued in this section aims to
represent the algorithm output in terms of desired or known quantities.
The argument has similarity to ones known in the Fourier theory of sampling.

\begin{mytheorem}
\label{thm:FiniteNresult_tilde_f}
Let $X_i$, $\hat F_n$, $Z$, $G$ and $\tilde f$ be given as
in section \ref{sec:CommonDefinitions}. % equations \eqref{eqn:G_definition} and \eqref{eqn:tilde_f_definition}. 
Then
\begin{align}
    \tilde f = \tfrac{T}{N} \cdot \sum_{z \in \{Z\}} p_z \cdot R_{N,T}((\cdot) - z) \label{eqn:tilde_f_is_periodic_smoothed_f_Z}
\end{align}    
where %$f_Z$ is the tempered distribution (in the sense of: ''generalized function'') representing the probability density of $Z$, and 
\begin{align}
    R_{N,T}(x) := \tfrac{1}{T} \cdot \frac{\sin(2 \pi \tfrac{2N-1}{2 T} \cdot x)}{\sin(2 \pi \tfrac{1}{2T} \cdot x)}  \label{eqn:R_N_T_definition}
\end{align}
with the expression on the right deemed continuously continued at the zeros of the denominator.
\end{mytheorem}

\begin{myremark}
The statement of the theorem is also written as 
\begin{align}
    \tilde f(x) = \tfrac{T}{N} \cdot \brackets{ f_Z \conv R_{N,T} }(x)  \label{eqn:ConvolutionalRepresentation_tilde_f_x}
\end{align}
where ''$*$'' signifies the convolution of two tempered distributions (i.e. in the sense of 
a generalized function, see appendix \ref{sec:FunctionalAnalyticBackground})
, and 
\begin{align}
    f_Z = \sum_{z \in \{Z\}} p_z \cdot \delta_z
\end{align}
is the tempered distribution 
representing the probability density belonging to $F_Z$. (The $\delta_z$ are delta distributions
with mass at $z$.)
\end{myremark}

\mbox{ }\\
Proof of the theorem: 

Utilizing the representation of $F_Z$, it is the Fourier transformation $G$ of $Z$ 
written as $G(\nu) = \E(\exp(- 2 \pi i \cdot Z \nu)) = \sum_{z \in \{Z\}} p_z \cdot \exp(- 2 \pi i \cdot z \nu)$.
Substituting this in the algorithm output \eqref{eqn:tilde_f_x_definition} yields
\begin{align}
    \tilde f(x) 
        & = \frac{1}{N} \sum_{k=-N+1}^{N-1} \sum_{z \in \{Z\}} p_z \cdot \exp(-2\pi i \cdot z k/T) \cdot \exp(2 \pi i \cdot x k/T) \\
        & = \frac{T}{N} \sum_{z \in \{Z\}} p_z \cdot \frac{1}{T} \sum_{k=-N+1}^{N-1} \exp(2 \pi i \cdot (x - z) k/T ) 
\end{align}
Identifying the finite geometric series in the sum over $k$ and applying an exponential
factor at it yields the above stated expression for $R_{N,T}$.
(The algorithm output is obtained by evaluating $\tilde f$ at the
places $i \tfrac{T}{N}$, $i=0 \dots N-1$.) \qed

The $R_{N,T}$ has properties which will allow to derive useful characteristics of the algorithm output.
It is $\lim_{x \converg 0+k T} R_{N,T}(x) = \tfrac{2N-1}{T}$ for $k \in \Whole$. 
It is $R_{N,T}$ periodic with period $T$, and $\int_0^T R_{N,T}(x) \d x = 1$. 
As $N \converg \infty$, have $R_{N,T} \converg \DiracComb_T$ in the space $\SchwartzSp'$ of tempered distributions
(see appendix \ref{sec:FunctionalAnalyticBackground}),
where $\DiracComb_T$ denotes the Dirac comb\footnote{Informal definition: 
$\DiracComb_T := \sum_{k=-\infty}^\infty \delta((\cdot)-kT)$. %Formal definition $\DiracComb_T(\phi) = \sum_{k=-\infty}^\infty \phi(kT)$.
} with interval $T$.
(The envelope of $R_{N,T}$, as $N$ increases, remains constant, but the increasing number
of zero crossings means that in intervals not overlapping with $T \cdot \Whole$, the positive and negative contributions cancel.)

\subsection{Motivating the proposed expression for the density estimate}
\label{sec:MotivatingTheDensityEstimateExpression}
According to Theorem \ref{thm:FiniteNresult_tilde_f}, the quantity $\tilde f$ can be regarded as proxy for a
density estimate, and (according to the remark subsequent to the common definitions)
itself fulfills a normalizing constraint $\sum_{i=0}^{N-1} \tilde f_i = 1$. 
Because this constraint implies that the $\tilde f_i$ scale to zero like $N^{-1}$ as 
$N \converg \infty$, a to-be-defined density estimate $\hat f$ must reasonably be 
of form $c \cdot N \cdot \tilde f$ in order to asymptotically fulfill the normalization.
The $c$ is found by noting that (or rather: aiming for)
\begin{align}
    1 \overset{!}{=} \int_{\kappa \cdot z_{min}}^{\kappa \cdot z_{min}+T} \hat f(x) \d x = \int_{\kappa \cdot z_{min}}^{\kappa \cdot z_{min}+T} c \cdot N \cdot \tilde f(x) \d x \approx c \cdot N \cdot \sum_{i=0}^{N-1} \tilde f_i \cdot \tfrac{T}{N},
\end{align}    
where the $\approx$ sign appeals to the numerical integration of the $\tilde f$ integral
using $N$ equidistant samples at $i\tfrac{T}{N}$. Thus reasonably $c=1/T$. 

Thus, if we here and henceforth set $I = [\kappa \cdot z_{min}, \kappa \cdot z_{max}]$, then the 
reasonable estimate for a smooth approximation to the density of $Z$ is
\begin{align}
    \hat f(x) := \tfrac{N}{T} \tilde f(x) \cdot 1_{I}(x) \label{eqn:final_Z_pdf_estimate}
\end{align}

\subsection{Asymptotics}
\label{sec:AsymptoticsBasedOnTildefRepresentation}
It shall now be proven, using the representation as periodic superposition
of the smoothed true density, that the integrated Algorithm 1 output converges
to the true cumulative distribution function as $N \converg \infty$. 
With ''integrated'' initially is meant the appropriate summation of the $\tilde f_i$ (regardable
as scaled density estimate according to the previous proposition).

\begin{mytheorem}
\label{thm:FiniteNbasedAsymptoticResult}
Let $X_i, a_j, Z, f_Z, F_Z, T_Z$ be given as in Theorem \ref{thm:FiniteNresult_tilde_f}, 
let $T$ be chosen $T > T_Z$, and $N$ be chosen.
Let $\tilde f$ be defined as in equation \eqref{eqn:tilde_f_x_definition} (equalling 
the expression \eqref{eqn:ConvolutionalRepresentation_tilde_f_x} below Theorem \ref{thm:FiniteNresult_tilde_f}). 
Call $I := [\kappa\cdot z_{min}, \kappa\cdot z_{max}]$ and $I_x := I \cap [\kappa\cdot z_{min}, x]$ (with $I_x = \emptyset$ for $x < \kappa\cdot z_{min}$), and set
\begin{align}
    \tilde F(x) := \int_{-\infty}^x \tfrac{N}{T} \cdot \tilde f(\xi) \cdot 1_{I}(\xi) \d \xi. \label{eqn:tilde_F_x_definition}
\end{align}
Let $A_d$ be the set of points of discontinuity of $F_Z$. Then it holds as $N \converg \infty$:
\begin{align}
    \text{i)} & \quad \tilde F(x) \converg F_Z(x) \quad \forall x \in \Real \backslash A_d \quad\quad\quad\quad\quad\quad\quad\quad\quad\quad\quad\quad\quad\quad\quad\quad\quad\quad\quad\quad\quad
\end{align}
\end{mytheorem}
\mbox{ }\\
Remark: The theorem (as the previous ones) considers the $X_i$ as given. Consequently also $F_Z$
is given as deterministic quantity. The convergence therefore is rightfully meant as a deterministic one.

\mbox{ }\\
Proof of (i): As before, $\SchwartzSp$ shall denote the (Schwartz) space of rapidly decreasing functions, 
and $\SchwartzSp'$ the associated space of tempered distributions. 
It was already stated that $R_{N,T} \converg \DiracComb_T$ in the distributional limit sense as $N \converg \infty$. 
Further, the convolution operation (among distributions) is continuous in the sense here needed\footnote{See results from appendix \ref{sec:FunctionalAnalyticBackground}.} It follows
\begin{align}
    \lim_{N \converg \infty} \tfrac{N}{T} \cdot \tilde f & = f_Z \conv \lim_{N \converg \infty} R_{N,T} = f_Z \conv \DiracComb_T.
\end{align}
Since for $T > T_Z$ the expression $f_Z \conv \DiracComb_T$ denotes the periodic
repetition of $f_Z$, it follows 
\begin{align}
    \tfrac{N}{T} \cdot \tilde f \cdot 1_{I} \converg (f_Z \conv \DiracComb_T) \cdot 1_{I} = f_Z. \label{eqn:tilde_f_x_masked_converges_to_f_Z}
\end{align}
Next, note that $F_Z(x) = \int_{-\infty}^x f_Z$ for all $x \not \in A_d$, where the 
integral of the distribution $f_Z$ over $[a,b]$ is defined via evaluation at a mollified $1_{[a,b]}$ (see 
appendix \ref{sec:FunctionalAnalyticBackground}). For any sequence of $\SchwartzSp$ functions 
$g_{n'}$ which converge (embedded in $\SchwartzSp'$) towards a distribution $G$, it is then $\int_{a}^b (g_{n'}) \converg \int_{a}^b G$
as ${n'} \converg \infty$. Consequently here also
\begin{align}
    \lim_{N\converg \infty} & \int_{-\infty}^x \tfrac{N}{T} \cdot \tilde f(\xi) \cdot 1_{I}(\xi) \d \xi = \lim_{N \converg \infty} \int_{I_x} \tfrac{N}{T} \cdot \tilde f  = \int_{I_x} \lim_{N\converg \infty} \tfrac{N}{T} \cdot \tilde f   \notag \\
        & = \int_{I_x} (f_Z \conv \DiracComb_T) = F_Z(x)
\end{align}
for all $x \in \Real \backslash A_d$.
This completes the proof of i). \qed

It is thus seen that the integral behind $\tilde F(\cdot)$ tends to the desired function $F_Z$
at all places of continuity of $F_Z$. One would now want to continue by proving a statement along 
\begin{align}
    \text{ii)} & \quad \tilde F_i - \tilde F(i T/N) \converg 0 \quad \forall i=0 \dots N-1 \quad\quad\quad\quad\quad\quad\quad\quad\quad\quad\quad\quad\quad\quad\quad\quad\quad
\end{align}
under the same assumptions as used for i), and where $\tilde F_i$ is the sum of the $\tilde f_i$ as in equation \eqref{eqn:tilde_F_i_definition}. 
The viewpoint behind this claim is that, in schematic words, $\tilde F_i$ is the sum of the $\tilde f_i = \tilde f(i\tfrac{T}{N})$ over suitable interval,
while $\tilde F(\cdot)$ is the integral of $\tilde f(\cdot)$ over the same interval. Thus $\tilde F_i$ roughly is 
the numerical evaluation (by equidistant sampling and rectangular rule) of the integral
behind $\tilde F(\cdot)$. A problem here occurs since as $N$ increases and the integration
partitioning becomes finer, also the integrand changes. It is then possible that the evaluations 
$\tilde f_i$ of $\tilde f$ all are errorenous, with the errors accumulating in the sum representing 
$\tilde F_i$. This will occur particularly if the possible
values of $Z$, and thus the placement of the $R_{N,T}(\cdot)$, exhibit a regularity.  
(Unrelated to the regularity in $f_Z$, the sampling of a 
smoothed $f_Z$ can be expected to produce errorenous results for certain choices 
of $N$ in a way that exhibits oscillation as $N$ varies.) 

To naturally avoid this, it is more suitable to sample $\tilde F(\cdot)$. With $x_0 := \kappa \cdot z_{min}$, have
\begin{align}
    \tilde F(x) & = \frac{1}{T} \cdot \sum_{k=-N+1}^{N-1} g_k \int_{x_0}^x \exp(2 \pi i \cdot \xi k \Delta \nu) \d \xi \notag \\
        & = \frac{x-x_0}{T}  + \frac{1}{T} \cdot \sum_{0 < \abs{k} \leq N-1} g_k \int_{x_0}^x \exp(2 \pi i \cdot \xi k \Delta \nu) \d \xi \notag \\
        & = \frac{x-x_0}{T} + \sum_{0 < \abs{k} \leq N-1} g_k \frac{\exp(2 \pi i \cdot x k \Delta \nu) - \exp(2 \pi i \cdot x_0 k \Delta \nu)}{2 \pi i \cdot k}.
\end{align}
Evaluated at $x = i T/N$, and using $x_0 = i_{min} T/N$, obtain 
\begin{align}
    \tilde F(i T/N) & = \frac{i-i_{min}}{N} + \sum_{0 < \abs{k} \leq N-1} g_k  \cdot \frac{e^{2 \pi i \cdot i k/N} - e^{2 \pi i \cdot i_{min} k/N}}{2 \pi i \cdot k} \notag \\
                & = \frac{i-i_{min}}{N} + \frac{1}{N} \cdot \sum_{0 < \abs{k} \leq N-1} g_k \cdot \frac{e^{2 \pi i \cdot i k/N} - e^{2 \pi i \cdot i_{min} k/N}}{2 \pi i \cdot k/N} \notag \\
                & = \frac{i-i_{min}}{N} + 2 \text{Re}\brackets{ \frac{1}{N} \cdot \sum_{k=1}^{N-1} g_k \cdot \frac{e^{2 \pi i \cdot i k/N}  - e^{2 \pi i \cdot i_{min} k/N}}{2 \pi i \cdot k/N} }.
\end{align}
This looks remarkably similar to the expression already available by the summation 
of the $\tilde f_i$. It is  
\begin{align}
    \tilde F_{i} & = \sum_{i'=i_{min}}^{i-1} \tilde f_{i} = \frac{1}{N} \sum_{k=-N+1}^{N-1} g_k \cdot \sum_{i'=0}^{i-i_{min}-1} \exp(2 \pi i \cdot \tfrac{(i'+i_{min}) k}{N})  \label{eqn:tilde_F2_i_definition} \\
        & = \frac{i-i_{min}}{N} + \frac{1}{N} \sum_{0 < \abs{k} \leq N-1} g_k \cdot \frac{e^{2 \pi i \cdot \tfrac{i k}{N}} - e^{2 \pi i \cdot \tfrac{i_{min} k}{N}}}{\exp(2\pi i \cdot k/N)-1} \\
        & = \frac{i-i_{min}}{N} + 2 \text{Re}\brackets{ \frac{1}{N} \sum_{k=1}^{N-1} g_k \cdot \frac{e^{2 \pi i \cdot \tfrac{i k}{N}} - e^{2 \pi i \cdot \tfrac{i_{min} k}{N}}}{\exp(2\pi i \cdot k/N)-1} }. \label{eqn:tilde_F_i_evaluation2}
\end{align}
One concludes that the desired values of $\tilde F(i T/N)$ can be obtained 
by computing the $\tilde f_i$ as in the previous way, \emph{but with the $g_k$ 
suitably modified by a factor}.

Concretely, denoting by $\dtilde F_i$ the expression constructed from $\tilde F_i$ 
by replacing $G(k \Delta \nu)$ in \eqref{eqn:tilde_F_i_definition} 
with $G(k \Delta) \cdot (\exp(2\pi i\cdot k/N)-1)/(2\pi i \cdot k/N)$, using the limit ''$k \converg 0$'' 
of the fraction at $k=0$, one obtains
\begin{align}
    \dtilde F_i = \tilde F(i T/N). \label{eqn:dtilde_F_result1}
\end{align} 
\begin{mycommenttext}
(With regards to the reasoning on the error convergence, this modification allowed 
to circumvent the part of the proof requiring stronger assumptions. Further, crucially, 
this modification does not impede the relevance of the previous theorem.)
\end{mycommenttext}
This is recorded in the lemma below.\vspace{1mm}

The modified version of the Algorithm 1 is henceforth referred to as ''Algorithm 2''.

\begin{mylemma}
\label{thm:FiniteNbasedAsymptoticResultAdapted}
Let $X_i, a_j, Z, G, f_Z, F_Z, T_Z$ be given as in Theorem \ref{thm:FiniteNbasedAsymptoticResult},
let $T$ be chosen as $T > T_Z$. Let $\tilde f$ be defined as in \eqref{eqn:tilde_f_x_definition}, and $\tilde F$ as in \eqref{eqn:tilde_F_x_definition}.
Let $\dtilde F_i$ be as defined before equation \eqref{eqn:dtilde_F_result1}. %as in equation \eqref{eqn:dtilde_F_definition}. 
Then
\begin{align}
    \text{ii')} & \quad \dtilde F_i - \tilde F(i T/N) = 0 \quad \forall i=i_{min} \dots i_{min}+N-1. \quad\quad\quad\quad\quad\quad\quad\quad\quad\quad\quad\quad\quad
\end{align}
\end{mylemma}
\mbox{ }\\
Proof: As argued before the statement of the lemma.\qed

Using i) of Theorem \ref{thm:FiniteNbasedAsymptoticResult} and the above ii'), one concludes 
that for all $i=0\dots N-1$ with $x = i \tfrac{T}{N} \not \in A_d$, it holds $\dtilde F_i \converg F_Z(i \tfrac{T}{N})$ as $N \converg \infty$.
Clearly, using $N$ sufficiently large, it will suffice to evaluate $\tilde F$
at those discrete places for which the theorem does provide the convergence assertion.
In practice and when not targetting specific assertions on the error, a value of $N = 1000$, in some 
applications $N=10000$, appears sufficient to compute for example quantiles of the distribution of $Z$. 
A derivation of a rigorous bound on the error at finite $N$ is found in the appendix \ref{sec:RateOfUniformConvergence}.

\section{Conclusion}
\label{sec:Conclusion}
The convergence of a variant of the algorithm in \cite{Pitschel2019DeterministicBootstrapping} (see 
modification derived in section \ref{sec:AsymptoticsBasedOnTildefRepresentation}) has been proven 
and a bound on the absolute value of error stated in dependence on the algorithm resolution $N$.
Future research could target further improvement of the convergence behaviour, in
particular by more effectively using the knowledge about the periodicity in the error component.
Moreover, details of the application of the algorithm to the areas mentioned in the introduction 
could be examined.

\appendix    
\section{Functional-analytic background for proof of Thm. \ref{thm:FiniteNbasedAsymptoticResult}}
\label{sec:FunctionalAnalyticBackground}
The objects in equation \eqref{eqn:ConvolutionalRepresentation_tilde_f_x} and in the proof of 
Theorem 2 part (i) are deemed elements of $\SchwartzSp'$, the
set of tempered distributions \cite{ReedSimon1972MathMethodsI}, i.e. continuous linear functionals $\SchwartzSp \converg \Complex$.
Here $\SchwartzSp$ is the space of rapidly decreasing $C^\infty$ functions $\Real \converg \Complex$, endowed with the
family of semi-norms $\norm{\phi}_{\alpha,\beta} := \max_\Real \abs{x^\beta D^\alpha \phi}$ (see \cite{ReedSimon1972MathMethodsI}).
Even though the in Theorem \ref{thm:FiniteNresult_tilde_f} mentioned $R_{N,T}(\cdot)$ is not 
in $\SchwartzSp$, regarding $R_{N,T}(\cdot)$ as tempered distribution allows
\begin{align}
    \FT[R_{N,T}](\phi) = (R_{N,T})(\FT[\phi]) = \int_\Real R_{N,T}(x) \FT[\phi](x) \d x, \quad \phi \in \SchwartzSp
\end{align}    
with converging integral on the right-hand side, so $\FT[R_{N,T}]$ is a well-defined
functional on $\SchwartzSp$, and is (because of continuity of $\FT$ in $\SchwartzSp$, see \cite{ReedSimon1972MathMethodsII}) continuous.
Thus $\FT[R_{N,T}] \in \SchwartzSp'$, and moreover its limit in $\SchwartzSp'$ is well-defined.
(In the following text, reference to ''tempered'' will be dropped, 
even though meant.)

For a distribution $F$ and a $g \in \SchwartzSp$, it is defined $(F \conv g)(\phi) := (F)((g \circ m) \conv \phi)$,
see \cite{ReedSimon1972MathMethodsII}. For $F \in \SchwartzSp'$, it will be $F \conv \delta_a$ be deemed defined via
approximating function sequence $\bar \delta_{a,1/s}$, knowing $\bar \delta_{a,1/s} \converg \delta_a$
as $s \converg \infty$. Then easily confirmed: $F \conv \delta_0 = F$, and $f(\cdot) \conv \delta_a = f((\cdot) + a)$,
and desirable properties (e.g. the convolution/product theorem) carry over via continuity.

More generally, for any distribution $G$ to which a sequence of $\SchwartzSp$ functions
$g_n$ converges (in $\SchwartzSp'$), can define $(F * G) := \lim_{n\converg\infty} (F * (g_n))$.
Then $G \mapsto (F \conv G)$ is continuous for those sequences, i.e. from $\lim_{n\converg\infty} (g_n) = G$
can conclude $\lim_{n\converg \infty} (F \conv (g_n)) = (F \conv \lim_{n\converg \infty} (g_n)) = (F \conv G)$.

The integral of a distribution over a bounded interval is defined 
as follows: for $a < b$, let $\bar J_{[a,b],\epsilon}$ be the function $1_{[a,b]}(\cdot) \conv \bar \delta_{0,\epsilon}$
with $\bar \delta_{0,\epsilon}$ as in proof of Theorem \ref{thm:FiniteNresult_tilde_f}.
It is $\bar J_{[a,b],\epsilon} \in \SchwartzSp$. Then set for $F \in \SchwartzSp'$, if the limit exists
and is independent of the shape of the mollifier choice $\bar \delta_{0,\epsilon}$,
\begin{align}
    \int_{a}^b F := \lim_{\epsilon \converg 0} (F)(\bar J_{[a,b],\epsilon}).
\end{align}
Clearly, for $F = (f)$ i.e. a distribution associated to an $L^1_{loc}$ function $f$,
it is $\int_{a}^b F$ existing and equal to the ordinary Lebesgue integral of $f$ over $[a,b]$.

\subsection{Convergence of $R_{N,T}$}
In section \ref{sec:FiniteNresults}, it was stated without proof that $R_{N,T} \converg \DiracComb_T$
as $N \converg \infty$. This is seen %when inspecting the intermediate steps in the proof there, or by explicitly considering the $R_{N,T}$, 
as follows: It was
$R_{N,T}(x) = T^{-1} \cdot \sum_{k=-N+1}^{N-1} \exp(2 \pi i \cdot k \Delta \nu x)$, 
therefore get (as distributional equation)
\begin{align}
    \lim_{N \converg \infty} & \FT^{-1}[\FT[R_{N,T}]] = \FT^{-1}[\lim_{N \converg \infty} \FT[R_{N,T}]] \\
            & = \FT^{-1}[\tfrac{1}{T} \cdot \lim_{N \converg \infty} \sum_{k=-N+1}^{N-1} \delta_{k \Delta \nu}] \\
            & = \FT^{-1}[\tfrac{1}{T} \cdot \DiracComb_{\Delta\nu}] = \FT^{-1}[\tfrac{1}{T} \cdot \DiracComb_{1/T}] = \DiracComb_{T},
\end{align}            
where the last equation uses a well-known result for the Dirac comb $\DiracComb_{T}$.
\footnote{The definition of the Fourier transformation $\FT$ on distributions 
underlying this equation is rooted in the definition of the Fourier transformation 
on the Schwartz space $\SchwartzSp$, using, for $v \in \SchwartzSp$, 
\begin{align}
    \FT[v](\nu) & := \int_\Real v(x) e^{-2 \pi i \cdot x \nu} \d x, \text{  and its inverse} \\
    \FT^{-1}[\hat v](x) & := \int_\Real \hat v(\nu) e^{2 \pi i \cdot x \nu} \d \nu.
\end{align}
}

\section{Rate of uniform convergence of $\tilde F(\cdot)$}
\label{sec:RateOfUniformConvergence}
\subsection{Overview}
The section gathers the elementary results leading up to 
the result
\begin{align}
    \abs{\tilde F(x) - F_Z(x)} = \ONotSym(N^{-1/2}) \text{ for each $x$ not a discontinuity place of $F_Z$}.
\end{align}    
It shall denote in this appendix $\multisetZ$ the \emph{multi}set of all possible values of $Z$ (each attained with probability $N_Z^{-1}$),
and $N_Z$ its cardinality. It is then $f_Z$ simply $f_Z = N_Z^{-1} \cdot \sum_{z \in \{Z\}} \delta_z$. 
Let $T > T_Z$, $\kappa = T/T_Z$, $I = [\kappa \cdot z_{min}, \kappa \cdot z_{max}]$ as before, 
and $x_0 = \kappa \cdot z_{min}$. By equation \eqref{eqn:tilde_F_x_definition}, it is 
\begin{align}
    \tilde F(x) = \int_{x_0}^x (f_Z \conv R_{N,T})(\xi) \d \xi,
\end{align}    
where the meaning of the integral as ordinary integral is justified
because the integrand is a function. Rewriting the integrand,
one obtains
\begin{align}
    \tilde F(x) = \int_{x_0}^x (\frac{1}{N_Z} \cdot \sum_{z \in \multisetZ} \delta_z \conv R_{N,T})(\xi) \d \xi = \frac{1}{N_Z} \cdot \int_{x_0}^x \sum_{z \in \multisetZ} R_{N,T}(\xi-z) \d \xi.
\end{align}    

On the other hand, for all $x \not \in \multisetZ$,
\begin{align}
    F_Z(x) = \int_{x_0}^x f_Z = \frac{1}{N_Z} \cdot \sum_{z \in \multisetZ} 1_{z \leq x}, % = \frac{1}{N_Z} \cdot \sum_{z \in \multisetZ} H(x-z).
\end{align}    
where the integral of the tempered distribution $f_Z$ is deemed defined 
as in appendix \ref{sec:FunctionalAnalyticBackground}.
For controlling the difference it therefore suffices to bound the 
terms 
\begin{align}
    \abs{\int_{x_0}^x R_{N,T}(\xi-z) \d \xi  \,\, - \,\, 1_{z \leq x}}
\end{align}    
for various $z \in \{Z\}$. To this end, the elementary integral
\begin{align}
    J_{N,T}(x) := \int_0^x R_{N,T}(\xi) \d \xi
\end{align}    
is further examined.

\subsection{Analysis of the integral of $R_{N,T}$}
It exhibits $R_{N,T}$ an oscillatory behaviour with an amplitude envelope 
which does not recede to zero (as $N \converg \infty$), no matter at which place $x$ 
this property is considered. (This is in contrast to the $\si(\cdot)$ function.)
However, since $R_{N,T}$ bounded by this envelope $T^{-1} \cdot (\sin(2\pi \xi/(2T)))^{-1}$
for all $\xi \not \in T \cdot \Whole$, its integral from $0$ to $x$ tends pointwise (for these $x$) to a function,
below seen to be $\tfrac{1}{2} + \gausslower{x/T}$. The following shows that in fact
the convergence is uniform on any closed interval within $[-T/2,T/2]$ not containing the zero.
The main tool employed for this is standard Fourier analysis.

Let $J_{N,T}$ be as above. Since $R_{N,T}$ is alternatively written as $R_{N,T}(x) = T^{-1} \cdot (1 + 2 \sum_{k=1}^{N-1} \cos(2 \pi k \cdot x/T))$,
it is $J_{N,T}$ equal to
\begin{align}
    J_{N,T}(x) = \frac{x}{T} + \frac{1}{\pi} \sum_{k=1}^{N-1} \frac{1}{k} \sin(2 \pi k \cdot x/T).
\end{align}    
In order to show $J_{N,T}(x) \converg \tfrac{1}{2} + \gausslower{x/T}$ pointwise, 
define the periodic odd function
\begin{align}
    h_T(x) := \tfrac{1}{2} - \frac{x}{T} + \gausslower{x/T} \text{ for $x \not \in T \cdot \Whole$}
\end{align}    
with $h_T(x)$ (arbitrarily) set equal to zero otherwise. It suffices then to show that
\begin{align}
    A_{N,T}(x) := \frac{1}{\pi} \sum_{k=1}^{N-1} \frac{1}{k} \sin(2 \pi k \cdot x/T) \,\converg \, h_T(x) \quad \text{ as } N \converg \infty.
\end{align}    
For this, one recognizes that $A_{N,T}$ is the truncated Fourier series expansion of $h_T$.
(Set $T=1$ and employ the orthonormal system $\{v_k\}$ with $v_k(x') := \sqrt{2} \cdot \sin(2 \pi k \cdot x')$
on the interval $[0,1]$.) Standard Fourier analytical techniques then immediately 
yield convergence in the $L_2$-norm on $[0,T]$, since $h_T$ is piecewise continuous.

For obtaining the uniform convergence and a non-asymptotic (numerical) bound, 
a result building on the special form of the summands in $A_{N,T}$
can be employed. (See for example Theorem 6.5 (Abel's test) in \cite{Walker1988FourierAnalysis},
repeated below for convenience.) 

\begin{mytheorem} 
\label{thm:AbelsTest}
(condensed from \cite{Walker1988FourierAnalysis}) Let $(f_n)_{n \in \Nat}$ be a sequence of 
real-valued functions defined on $[a,b]$, fulfilling that the absolute value of the partial sums of $\sum_{n=1}^\infty f_n(x)$
are uniformly bounded, say by constant $M \in \Real$. Let $(a_n)_{n \in \Nat} \subset \Real$ be a decreasing sequence
converging to zero. Then $\sum_{n=1}^N a_n f_n(x)$ converges uniformly on $[a,b]$ as $N \converg \infty$, and
the residual sums fulfill for all $N \in \Nat$
\begin{align}
    \abs{ \sum_{n=N+1}^\infty a_n f_n(x) } \leq 2 a_{N+1} M. \label{eqn:AbelResidualSumBound}
\end{align}    
uniformly in $x$. (The proof of the theorem is based on Abel's lemma.) \qed  
\end{mytheorem}

\mbox{ }\\
The application of the theorem to the present case allows to derive a bound on the error
caused by truncation of the infinite series, yielding the following result.

\begin{mylemma}
\label{thm:UniformConvergenceOfJNToutsideZero}
Let $A_{N,T}$ be defined as above. Let $I_2$ be a closed interval $\subset \Real$ with
$\{\abs{x-j T}, x \in I_2, j \in \Whole\} \geq T \cdot \epsilon^* > 0$, i.e. $I_2$ is bounded
away from $T \cdot \Whole$ by at least $T \cdot \epsilon^*$. Then $A_{N,T}(x) \converg h_T(x)$ uniformly
on $I_2$, and for all $x \in I_2$
\begin{align}
    \abs{A_{N,T}^\infty(x)} := \abs{\frac{1}{\pi} \sum_{k=N}^\infty \frac{1}{k} \sin(2 \pi k \cdot x/T) } \leq \frac{1}{\pi N} \cdot \frac{1}{\epsilon^*}\,\,.  \label{eqn:SinNxResisualBound} %2019-06-04:eqn25
\end{align}    
\end{mylemma}
\mbox{ }\\
Proof: Similarly as in the proof of Theorem \ref{thm:FiniteNresult_tilde_f}, %a bound for $Q_1$
an expression for $Q_1 := \sum_{k=1}^{N-1} \sin(2\pi k \cdot x/T)$ can be obtained.
With $x' = x/T$ and $0 < x' < 1/2$ 
\begin{align}
    Q_1 & = \text{Im}\brackets{\sum_{k=1}^{N-1} \exp(2\pi i \cdot k \cdot x') } = \text{Im}\brackets{ \frac{e^{2 \pi i \cdot N x'} - e^{2 \pi i \cdot x'} }{ e^{2 \pi i \cdot x'} - 1 } } \\
        & = \frac{1}{2 \sin(2 \pi \cdot x'/2)} \cdot \text{Re}\brackets{e^{2 \pi i \cdot (N-\tfrac{1}{2}) x'} - e^{2 \pi i \cdot \tfrac{1}{2} x'}} \leq \frac{1}{\sin(2 \pi \cdot x'/2)}.
\end{align}    
With $\sin(u) \geq \tfrac{2}{\pi} u$ for $u \in [0,\pi/2]$ get
\begin{align}
    \abs{Q_1} \leq \frac{1}{(2/\pi) \cdot (2 \pi x'/2)} = \frac{1}{2 x'}
\end{align}    
(This bound is extended to $x' \in (0,1)$, by observing the symmetry of $\sin(2 \pi \cdot x'/2)$, 
yielding $Q_1 \leq 1 / (2 \min(x,1-x))$, and analogously to $x \in \Real \backslash (T \cdot \Whole)$.)
For $x \in I_2$, it follows
\begin{align}
    \abs{Q_1} \leq \frac{1}{2 \epsilon^*} =: M \quad \text{ for all } x \in I_2. \label{eqn:SinNxResidualBound}
\end{align}    
Applying theorem \ref{thm:AbelsTest} yields the result.\qed

\begin{myremark}
\label{rem:UniformConvergenceOfJNToutsideZero_reformulation}
Since the above arguments found that 
\begin{align}
    \frac{1}{2} + \gausslower{x/T} = \frac{x}{T} + \frac{1}{\pi} \sum_{k=1}^\infty \frac{1}{k} \sin(2 \pi k \cdot x/T),
\end{align}
the previous result becomes applicable on 
\begin{align}
    \abs{ J_{N,T}(x) - (\frac{1}{2} + \gausslower{x/T}) } = \abs{ -\frac{1}{\pi} \sum_{k=N}^\infty \frac{1}{k} \sin(2 \pi k \cdot x/T) }.
\end{align}
The left-hand side in turn equals $\abs{ J_{N,T}(x) - (-\frac{1}{2} + 1_{0 \leq x/T}) }$ for $x \in (-T,T)$. \qed
\end{myremark}

In the above, the factor $T$ appearing next to $\epsilon^*$ reflects
the fact that as we arbitrarily scale $Z$ (and let $\kappa$ be constant), 
the minimum distance of $I_2$ which needs to be kept from $T \cdot \Whole$
must vary linearly in $T$ in order to leave the bound for $\abs{Q_1}$
invariant. When employing the lemma, $\epsilon^*$ will be chosen (as usual) in dependence of $N$ 
such that the contribution to the total error from $J_{N,T}$ which turn
out to be evaluated in the $[-\epsilon^*, \epsilon^*]$ interval just balances 
with the contribution from the bound \eqref{eqn:SinNxResidualBound}. I.e., $\epsilon^*$
may not be choosen too small.

In the following theorem, the ''max'' expression appearing 
equals $\max_{z_0} P^Z([z_0 - T \epsilon^*, z_0 + T \epsilon^*])$, 
where $P^Z$ is the measure of $Z$ derived from $\hat F_n$.

\begin{mytheorem}
\label{thm:UniformConvergenceTildeFwithRateAndConstants}
Let $X_i, Z, F_Z, T_Z$ be as in section \ref{sec:CommonDefinitions} and $\tilde F$ be 
as in equation \eqref{eqn:tilde_F_x_definition}.
Let $T > T_Z$ and $\kappa = T/T_Z$, and $x_0$ chosen with $\kappa \cdot z_{min} \leq x_0 < z_{min}$. Let $M_2 \in \Real$ with
\begin{align}
    \max_{z_0 \in \Real} \frac{1}{N_Z} \!\! \sum_{\substack{z \in \multisetZ \\ \abs{z-z_0} \leq T \cdot \epsilon^*}} \!\!\! 1 \,\,\, & \leq  \,\, M_2 \cdot \epsilon^*. \label{eqn:M2DensityAssumption}
\end{align}
Then %for all $x \in I \backslash A_d = [\kappa \cdot z_{min}, \kappa \cdot z_{max}] \backslash A_d$,
\begin{align}
    \abs{\tilde F(x) - F_Z(x)} \leq 2 \cdot \sqrt{\frac{1}{2 \pi} M_2} \cdot N^{-1/2}
\end{align}
for each $x \in I \backslash A_d$, where $A_d$ is the set of places
of discontinuity of $F_Z$. 
\end{mytheorem}
\mbox{ }\\
Proof: Let $x \in I \backslash A_d$. It is, with $f_Z = \tfrac{1}{N_Z} \sum_{z \in \multisetZ} \delta_z$,
\begin{align}
    \abs{\tilde F(x) - F_Z(x)} & = \abs{\int_{x_0}^x \brackets{ \tfrac{1}{N_Z} \sum_{z \in \multisetZ} \delta_z \conv R_{N,T}}(\xi) \d \xi - \tfrac{1}{N_Z} \sum_{z \in \multisetZ} 1_{z \leq x} } \\
        & = \tfrac{1}{N_Z} \abs{\int_{x_0}^x \sum_{z \in \multisetZ} R_{N,T}(\xi-z) \d \xi - \sum_{z \in \multisetZ} 1_{z \leq x} } \\
        & = \tfrac{1}{N_Z} \abs{\sum_{z \in \multisetZ} (J_{N,T}(x-z) - J_{N,T}(x_0-z)) - \sum_{z \in \multisetZ} 1_{z \leq x} } \\
        & \leq \tfrac{1}{N_Z} \sum_{z \in \multisetZ} \abs{ (J_{N,T}(x-z) - J_{N,T}(x_0-z)) - 1_{z \leq x} } 
\end{align}
Since $x_0 < \min Z$ and $\max Z - x_0 < \kappa \cdot z_{max} - x_0 \leq T$, it is $x_0 - z < 0$ and 
$x_0 - z > x_0 - \max Z > - T$, i.e. $x_0 - z \in (-T,0)$ for all $z \in \multisetZ$.
Therefore the lower evaluation of $J_{N,T}$ tends to $-\tfrac{1}{2}$ as $N \converg \infty$.
With $x \not \in A_d$, it is $x \neq z$. Via $x \in I$, it is $x-z \in (-T,T)$.
For $x > z$, it tends $J_{N,T}(x-z)$ towards $\tfrac{1}{2}$. With this
\begin{align}
   \abs{\tilde F(x) - F_Z(x)} & \leq \tfrac{1}{N_Z} \sum_{z \in \multisetZ} \abs{ - J_{N,T}(x_0-z) - \tfrac{1}{2} } + \abs{ J_{N,T}(x-z) + \tfrac{1}{2} - 1_{z \leq x} } \\
        & \leq \tfrac{1}{N_Z} \sum_{z \in \multisetZ} \abs{ - J_{N,T}(x_0-z) - \tfrac{1}{2} } + \abs{ \frac{1}{\pi} \sum_{k=N}^\infty \frac{1}{k} \sin(2 \pi k \cdot (x-z)/T)) } \label{eqn:UniformConvergenceProof_eqnA}
\end{align}
where the remark \ref{rem:UniformConvergenceOfJNToutsideZero_reformulation} was used 
observing $x - z \in (-T,T)$.

Examine the second term first. For summands with $z$ close to $x$, Lemma 
\ref{thm:UniformConvergenceOfJNToutsideZero} will not be applicable. Thus let 
$\epsilon^* > 0$ and represent $\multisetZ$ as union of the multisets 
$A_1(x) := \{z \in \multisetZ, \abs{x-z} \leq T \cdot \epsilon^*\}$ and 
$A_2(x) := \{z \in \multisetZ, \abs{x-z} > T \cdot \epsilon^*\}$.
Then using the assumption \eqref{eqn:M2DensityAssumption} in the summation over $A_1$ yields
\begin{align}
        N_Z^{-1} & \sum_{z \in \multisetZ} \abs{ \frac{1}{\pi} \sum_{k=N}^\infty \frac{1}{k} \sin(2 \pi k \cdot (x-z)/T)) } \\
        & \leq N_Z^{-1} \cdot \sum_{z \in A_1(x)} \abs{A_{N,T}^\infty(x-z)} + N_Z^{-1} \cdot \sum_{z \in A_2(x)} \abs{A_{N,T}^\infty(x-z)} \\
        & \leq \frac{1}{2} \cdot M_2 \cdot \epsilon^* + \frac{1}{\pi N} \cdot \frac{1}{\epsilon^*} =: E_2.
\end{align}
Herein, it was used that $\abs{J_{N,T}(x-z) + \tfrac{1}{2} - 1_{z \leq x}}$ is bounded by $\tfrac{1}{2}$ uniformly in $N$ and $x$.
Similarly, the first term in eqn. \eqref{eqn:UniformConvergenceProof_eqnA} is bounded, yielding the same 
bound. \footnote{In fact, if $\kappa$ is sufficiently large or $\epsilon^*$ sufficiently small
such that $T \cdot \epsilon^* < (\kappa-1) z_{min}$, then $A_1(x_0)$ will be empty, 
which effects below a reduction in the ultimate bound by factor $1/\sqrt{2}$.} In total
\begin{align}
    \abs{\tilde F(x) - F_Z(x)} & \leq 2 \cdot \brackets{ \frac{1}{2} \cdot M_2 \cdot \epsilon^* + \frac{1}{\pi N} \cdot \frac{1}{\epsilon^*} } =: E.
\end{align}
It is easy to determine that the optimum choice for $\epsilon^*$ minimizing
$E$ yields 
\begin{align}
    E_{\min} = 2 \cdot \sqrt{\frac{1}{2} M_2 \cdot \frac{1}{\pi}} \cdot N^{-1/2},
\end{align}
which proves the theorem. \qed

The overall result obtained therefore is: For places $x$ in the discontinuity set $A_d$ of $F_Z$,
no statement is made. For places $x$ outside this set, Theorem 
\ref{thm:UniformConvergenceTildeFwithRateAndConstants} states a bound on the pointwise error
which uniformly holds over $I \backslash A_d$. 

\mbox{ }\\
Plausible values for $M_2$ can be only estimated crudely before having obtained an
estimate of $f_Z$. The following remark provides a hint to possible occuring values in 
the non-degenerate case (usually resulting from $X$ having a continuous distribution)
by deriving an $M_2$ for the normal distribution. In cases where the sample $(X_i)_{i=1...n}$ 
exhibits regularities (e.g. $X_i \in \Whole$ for all $i$) or even is degenerate ($X_i = \const$),
the $M_2$ will have excessive values.

\begin{myremark}
Let $U \sim \NVert(0, \sigma^2)$ and $T = 5 \sigma$. Let $0 < \epsilon^* < 5$. Then
\begin{align}
    \max_{z_0} P(\abs{Z' - z_0} \leq T \cdot \epsilon^*) & = \max_{z_0}  \frac{1}{\sqrt{2 \pi} \sigma} \int_{z_0-T \epsilon^*}^{z_0+T \epsilon^*} e^{-z^2/(2 \sigma^2)} \d z \\
            & \leq \frac{2 \epsilon^* \cdot T}{\sqrt{2 \pi} \sigma} = \frac{10}{\sqrt{2 \pi}} \epsilon^*,
\end{align}
i.e. for the random variable $U$ the appropriate value for $M_2$ would be $10/\sqrt{2 \pi}$. \qed
\end{myremark}

\kommentar{
\mbox{ }\\
In the following remark, it is 
$M_2(\tilde I)$ supposed as a bound on $\max_{z_0 \in \tilde I} P^Z([z_0 - T \epsilon^*, z_0 + T \epsilon^*])$,
i.e. extends the previously used $M_2$.

\begin{myremark}
If, instead of on all of the interval $I$, one ask for the
error bound in a region $\tilde I \subset I$ belonging to rare (i.e. usually: extremal) values of $Z$,
the corresponding $M_2$ restricted to that region, denoted as $M_2(\tilde I)$, will 
be substantially lower than $M_2(I)$, for example only a tenth of it. Such case will arise
for example when seeking quantiles for confidence intervals $\alpha$ with small $\alpha$.
As rule of thumb, the suitable $M_2$ at any point $x \in I \backslash A_d$ 
will be approximately proportional to the computed smoothed density at this point.
\end{myremark}
}

\bibliographystyle{alpha}
\bibliography{dbt_bib_man}{}

\end{document}